\documentclass[12pt]{amsart}
\usepackage[margin=1.25in]{geometry}

\usepackage{amsmath, amscd, amsthm}
\usepackage[pdftex]{color}
\usepackage{graphicx}
\usepackage[all]{xy}

\usepackage{hyperref}

\newtheorem{theorem}{Theorem}
\newtheorem{lemma}[theorem]{Lemma}
\newtheorem{proposition}[theorem]{Proposition}
\newtheorem{corollary}[theorem]{Corollary}

\newtheorem{definition}[theorem]{Definition}

\theoremstyle{definition}
\newtheorem{remark}[theorem]{Remark}

\title{Generalizations of Douady's magic formula}

\author{Adam Epstein}
\thanks{University of Warwick, UK, \url{A.L.Epstein@warwick.ac.uk}}

\author{Giulio Tiozzo}
\thanks{University of Toronto, Canada, \url{tiozzo@math.utoronto.ca}}

\begin{document}

\begin{abstract}
We generalize a combinatorial formula of Douady from the main cardioid to other hyperbolic components $H$ of the Mandelbrot set, 
constructing an explicit piecewise linear map which sends the set of angles of external rays landing on $H$ 
to the set of angles of external rays landing on the real axis.
\end{abstract}

\maketitle

In the study of the dynamics of the quadratic family $f_c(z) := z^2 + c$, $c \in \mathbb{C}$, a central object of interest is the \emph{Mandelbrot set} 
$$\mathcal{M} := \{ c \in \mathbb{C} \ : \ (f_c^{(n)}(0))_{n \geq 0} \textup{ is bounded} \}.$$

The Mandelbrot set contains two notable smooth curves: namely, the intersection with the real axis $\mathcal{M} \cap \mathbb{R} = [-2, 1/4]$, and the \emph{main cardioid} of $\mathcal{M}$, which is the set of parameters $c$ for which $f_c$ has an attracting or indifferent fixed point (of course, this is smooth except at the cusp $c = 1/4$). 

Let $\mathcal{R}$ denote the set of angles of external rays landing on the real axis, and $\Xi$ the set of angles of external rays landing on the main cardioid. The following ``magic formula" is due to Douady (for a proof see \cite{Bl}): 

\begin{theorem}[Douady] \label{T:Douady}
The map 
$$T(\theta) := \left\{ \begin{array}{ll} 
\frac{1}{2} + \frac{\theta}{4} & \textup{if }0 \leq \theta < \frac{1}{2} \\
\frac{1}{4} + \frac{\theta}{4} & \textup{if }\frac{1}{2} < \theta \leq 1
\end{array}\right.$$ 
sends $\Xi$ into $\mathcal{R}$.
\end{theorem}

In this note, we prove the following generalization of Theorem \ref{T:Douady} to other hyperbolic components. Let $D(\theta) := 2 \theta \mod 1$
denote the doubling map; moreover, given a finite binary word $S$ and a real number $\theta$, by $S \cdot \theta$ we denote the real number whose binary expansion is the concatenation of $S$ and the expansion of $\theta$. 

\begin{theorem} \label{T:main}
Let $H$ be a hyperbolic component in the upper half plane which belongs to a vein $V$ of the Mandelbrot set, not in the $1/2$-limb. 
Let $\bold{A}_H$, $\bold{B}_H$ be the binary expansions of the root of $H$, with $\bold{A}_H < \bold{B}_H$, and let $\delta_V$ be the complexity of $V$. 
Then the map 
$$\Phi_H(\theta) := D^{\delta_V}(\bold{B}_H \bold{A}_H \cdot \theta)$$
sends the set of external angles of rays landing on the upper part of $H$ into the set $\mathcal{R}$ of angles of rays landing on the real axis.
\end{theorem}

\begin{figure} 
\begin{center}
\fbox{\includegraphics[width = 0.7 \textwidth]{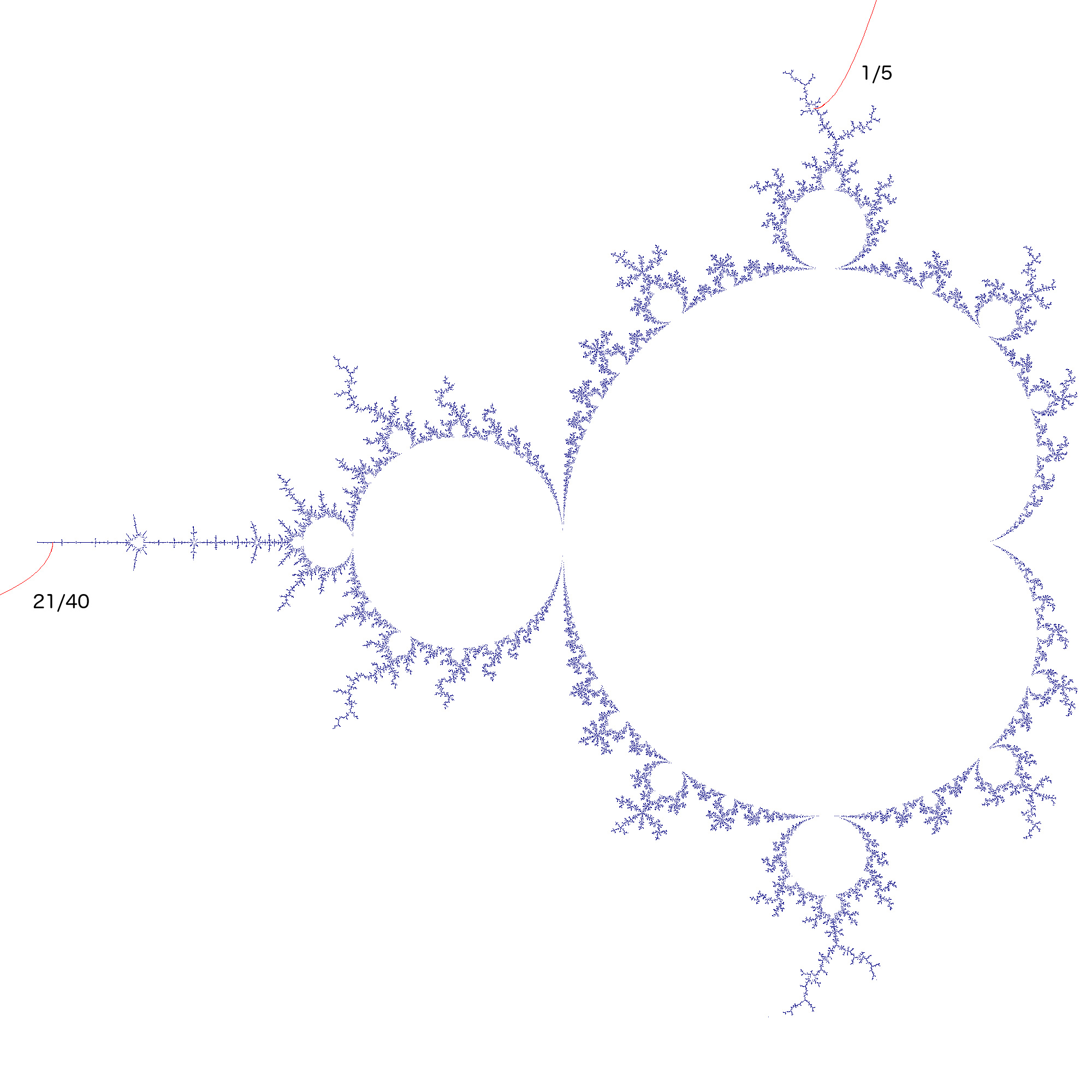}}
\end{center}
\caption{The magic formula for the kokopelli component. The angle $\theta = 1/5$ which lands at the root of its component 
is sent to $\Phi_H(\theta) = 21/40$ which lands on the real axis.}
\label{F:one}
\end{figure}

To illustrate such a formula, let us consider the ``kokopelli" component of period $4$, which has root of angles $\theta_1 = 3/15 = .\overline{0011}$, 
$\theta_2 = 4/15 = .\overline{0100}$, thus  $\bold{A}_H = 0011$, $\bold{B}_H = 0100$. This component lies on the principal vein with tip $\theta = 1/4$, hence its complexity is $\delta_V = 1$.
Thus, 
$$\Phi_H(\theta) := 1000011 \cdot \theta =  \frac{67}{128} + \frac{\theta}{128}.$$
Let us note that if $H$ is the main cardioid, then $\delta_V = 0$, $\bold{A}_H = 0$, $\bold{B}_H = 1$, hence we recover Douady's original formula from Theorem \ref{T:Douady}. By symmetry, an analogous formula to Theorem \ref{T:main} holds for hyperbolic components in the lower half plane.
 
\smallskip
Note that the map $\Phi_H$ is very far from being surjective: indeed, the set of angles of rays landing on a hyperbolic component 
has Hausdorff dimension zero, while the set $\mathcal{R}$ has dimension $1$.

We will prove our formula (Theorem \ref{T:main}) in two parts: first (Section \ref{S:BC}),  we will construct a map which sends the hyperbolic component into 
the real slice of the corresponding small Mandelbrot set; then (Section \ref{S:veins}), we will map the vein to which this small Mandelbrot set belongs to the real vein. This will be done by developing some notion of combinatorial veins and combinatorial Hubbard trees, 
which may be of independent interest.
The final formula will be the composition of the two formulas. 
Then, in Section \ref{S:renorm} we will prove that the rays given by the image of the formula actually land, by showing they are not renormalizable. 
Finally, in Section \ref{S:alter} we will prove an alternate formula, which, however does not depend on the vein structure of $\mathcal{M}$. 

During the preparation of this paper we have been informed of the generalization of Douady's formula by Bl\'e and Cabrera \cite{BC}. 
Let us remark that our formula is quite different, as the one in \cite{BC} does not map into the real axis but rather into the tuned copies 
of the real axis inside small Mandelbrot sets. In fact, we will describe their formula and how it relates to ours in Section \ref{S:BC}.

\begin{figure}
\begin{center}
\includegraphics[width = 0.7  \textwidth]{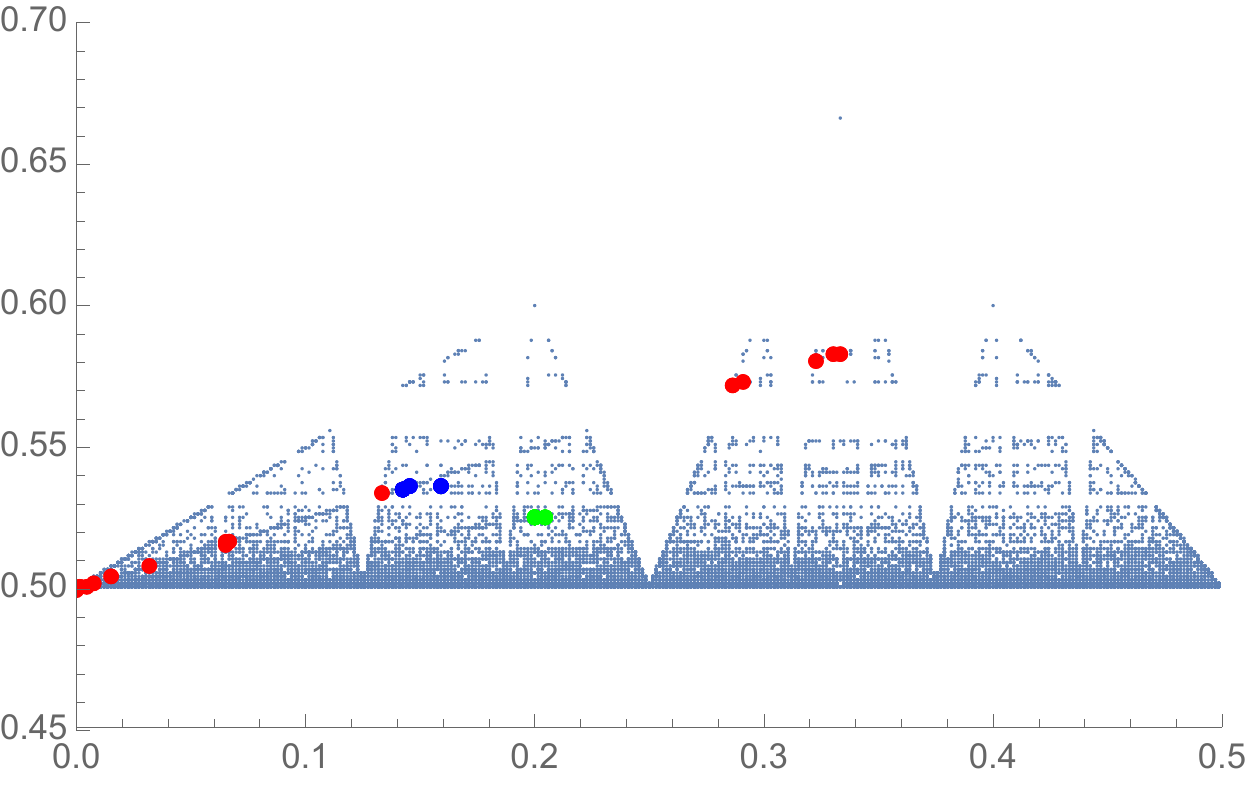}
\caption{The graph of the function $\Psi(x)$ 
discussed in Remark \ref{R:psi}. The red dots correspond to points on the graph of Douady's magic formula function $T = \varphi_H$ for the main cardioid, the (thick) blue dots to the image of $\varphi_H$ for the rabbit component, and the green dots to the kokopelli component.}
 \label{F:two}
\end{center}
\end{figure}

\subsection*{Acknowledgements} 
We thank S. Koch and C. T. McMullen for useful conversations. 
G.T. is partially supported by NSERC and the Alfred P. Sloan Foundation. 

\section{Combinatorics of external rays}

\subsection{Real rays and the original formula}

Consider the Riemann map $\Phi_M : \mathbb{C}~\setminus~\overline{\mathbb{D}} \to \mathbb{C} \setminus \mathcal{M}$, 
which is unique if we normalize it so that $\Phi_M(\infty) = \infty$, $\Phi_M'(\infty) = 1$.
For each $\theta \in \mathbb{R}/\mathbb{Z}$, we define the \emph{external ray} of angle $\theta$ as 
$$R_M(\theta) := \{ \Phi_M(r e^{2 \pi i \theta}) \ : \ r > 1 \}.$$
The ray $R_M(\theta)$ \emph{lands} if there exists $\gamma(\theta) := \lim_{r \to 1^+} \Phi_M(r e^{2 \pi i \theta})$. 
It is conjectured that the ray $R_M(\theta)$ lands for any $\theta \in \mathbb{R}/\mathbb{Z}$; 
to circumvent this issue, we define as $\widehat{R}_M(\theta)$ the impression of $R_M(\theta)$, 
and the set of \emph{real angles} as 
$$\widehat{\mathcal{R}} := \{ \theta \in \mathbb{R}/\mathbb{Z} \ : \ \widehat{R}_M(\theta) \cap \mathbb{R} \neq \emptyset \}.$$
Conjecturally, this is exactly the set of angles for which the corresponding ray lands on the real axis. For further details, see e.g. \cite{Za}.

The following lemma gives a characterization of $\widehat{\mathcal{R}}$ in terms of the dynamics of the doubling map.

\begin{lemma}[\cite{Do}] \label{L:real}
The set of real angles 
equals 
$$\widehat{\mathcal{R}} := \{ \theta \in \mathbb{R}/\mathbb{Z} \ : \  |D^n(\theta) - 1/2| \geq |\theta - 1/2| \quad \forall n \in \mathbb{N} \}.$$
\end{lemma}

\begin{proof}[Proof of Theorem \ref{T:Douady}]
By symmetry, we can assume $0 < \theta < \frac{1}{2}$.
Since the external ray of angle $\theta$ lands on the main cardioid, then the forward orbit $P := (D^n (\theta))_{n \in \mathbb{N}}$
of the angle $\theta$ does not intersect the half-circle $I = (\frac{\theta}{2} + \frac{1}{2}, \frac{\theta}{2})$ (the one which contains $0$) (see \cite{BS}).
The preimage of $I$ by the doubling map is the union of two intervals $(-\frac{\theta}{4}, \frac{\theta}{4}) \cup (\frac{1}{2} - \frac{\theta}{4}, \frac{1}{2} + \frac{\theta}{4})$. In particular, the forward orbit of $\theta$ does not intersect $J = (\frac{1}{2} - \frac{\theta}{4}, \frac{1}{2} + \frac{\theta}{4})$. Thus, if we look at the orbit of $\theta' = T(\theta)$, we have $2 \theta' = \frac{\theta}{2} \notin J$ and 
$D^n(\theta') \in P$ for any $n \geq 2$, hence the forward orbit of $\theta'$ does not intersect $J$, hence $\theta'$ is a real angle 
by Lemma \ref{L:real}. 
\end{proof}

\begin{remark} \label{R:psi}
By Lemma \ref{L:real}, one recognizes that a way to map any angle $\theta \in \mathbb{R}/\mathbb{Z}$ to the 
set of real angles is to consider the function 
$$\Psi(x) := \frac{1}{2} + \inf_{k \geq 0} \left| D^k(x) - \frac{1}{2} \right|,$$
which indeed satisfies $\Psi(\mathbb{R}/\mathbb{Z}) \subseteq \widehat{\mathcal{R}}$.
As you can see from Figure \ref{F:two}, this function is discontinuous, while its graph ``contains" the graphs of 
all the magic formula functions $\varphi_H$ given by Theorem \ref{T:main} (which are indeed continuous) for all components $H$.
\end{remark}

\subsection{Tuning and the Bl\'e-Cabrera magic formula} \label{S:BC}

Given any hyperbolic component $H$ in the Mandelbrot set, let us recall that there is a \emph{tuning map} which sends
the main cardioid to the hyperbolic component $H$, and the Mandelbrot set to a small copy of itself which contains $H$.

In order to define the map precisely, let $\bold{a}_H = .\overline{a_1\dots a_k} < \bold{b}_H =.\overline{b_1 \dots b_k}$ denote the two external angles of rays landing at the root 
of $H$, and denote $\bold{A}_H = a_1\dots a_k$ and $\bold{B}_H = b_1 \dots b_k$ the two corresponding finite binary words.
Moreover, denote $\bold{a}'_H = .\overline{a_1\dots a_k b_1 \dots b_k}$ and $\bold{b}'_H = .\overline{b_1 \dots b_k a_1\dots a_k }$.

The tuning map $\mathfrak{T}_H$ on the set of external angles is now defined as follows.
If $\theta = . \epsilon_1 \epsilon_2 \dots$ is the binary expansion of $\theta$, then the angle $\mathfrak{T}_H(\theta)$ has binary expansion
$$\mathfrak{T}_H(\theta) = .A_{\epsilon_1} A_{\epsilon_2} \dots$$
where $A_0 = \bold{A}_H$, $A_1 = \bold{B}_H$.
Then, the set $\Xi_H := \mathfrak{T}_H(\Xi)$ is the set of rays landing on the boundary of $H$.

If $\theta \in \mathbb{T} = \mathbb{R}/\mathbb{Z}$ has infinite binary expansion (i.e. it is not a dyadic number) and $S = s_1\dots s_n$ is a finite word on the alphabet $\{0,1\}$, 
we denote by $S \cdot \theta$ the element of $\mathbb{T}$ 
$$S \cdot \theta := \sum_{k=1}^n s_k 2^{-k} + 2^{-n} \theta$$
i.e. the point whose binary expansion is the concatenation of $S$ and the binary expansion of $\theta$.
Recall that tuning behaves well with respect to concatenation, namely for any $S, H$ and $\theta$ we have 
$$\mathfrak{T}_H(S \cdot \theta) = \mathfrak{T}_H(S) \cdot \mathfrak{T}_H(\theta).$$

We define the \emph{small real vein} associated to $H$, and denote it as $\mathbb{R}_H$, to be the real vein of the small Mandelbrot set associated to $H$.
Let us denote as $\mathcal{R}_H$ the set of external angles of rays landing on $\mathbb{R}_H$, and as $\widehat{\mathcal{R}}_H$ the set of external 
angles of rays whose impression intersects $\mathbb{R}_H$. Clearly, $\mathbb{R}_{H_0} = \mathbb{R}$ if $H_0$ is the main cardioid, and $\mathcal{R}_{H_0} = \mathcal{R}$.

\begin{proposition}[Bl\'e-Cabrera \cite{BC}] \label{P:BC-formula}
The map $T_H := \mathfrak{T}_H \circ T \circ \mathfrak{T}_H^{-1}$ maps $\Xi_H$ into the small real vein $\widehat{\mathcal{R}}_H$.
Moreover, this map (restricted to $\Xi_H$) is piecewise affine: in fact, it can be written in terms of binary expansions as 
$$T_H(\theta) = \bold{B}_H \bold{A}_H  \cdot \theta$$
if $\theta \in (\bold{a}_H, \bold{a}'_H)$, and 
$$T_H(\theta) =  \bold{A}_H  \bold{B}_H \cdot \theta$$
if $\theta \in (\bold{b}'_H, \bold{b}_H)$.
\end{proposition}

\begin{proof}
By construction, the tuning map $\mathfrak{T}_H : \Xi \setminus \{ 0\} \to \Xi_H \setminus \{ \bold{a}_H, \bold{b}_H \}$ is a bijection, so the first statement follows by looking at the diagram: 
$$\xymatrix{\Xi \ar[r]^T \ar[d]^{\mathfrak{T}_H} & \widehat{\mathcal{R}} \ar[d]^{\mathfrak{T}_H} \\
\Xi_H \ar[r]^{T_H} & \widehat{\mathcal{R}}_H. }$$
Moreover, since $T(\theta) = 10 \cdot \theta$ for $\theta <1/2$ and tuning behaves nicely with respect to concatenation, we have 
$$T_H(\theta) = \mathfrak{T}_H ( 10 \cdot \mathfrak{T}_H^{-1}(\theta) ) = \mathfrak{T}_H(10) \cdot \theta $$
which proves the second claim.
\end{proof}

\subsection{Veins, pseudocenters and complexity} \label{S:veins}

Given a dyadic number $\theta \in S^1 = \mathbb{R}/\mathbb{Z}$, we define its \emph{complexity} as
$$\Vert \theta \Vert := \min \{ k \geq 0 \ : \ D^k(\theta) = 0 \mod 1 \}.$$
Of course, if $\theta = \frac{p}{2^q}$ with $p$ odd, then $\Vert \theta \Vert = q$. 
Given an interval $(\theta^-,\theta^+)$ with $\theta^- < \theta^+$, we define its \emph{pseudocenter} $\theta_0$
to be the dyadic rational of lowest complexity inside the interval $(\theta^-, \theta^+)$. 

\medskip
A pair of elements $(\theta^-, \theta^+)$ in $\mathbb{T}$ is a \emph{ray pair} if the two external rays of angle $\theta^-$ and 
$\theta^+$ combinatorially land at the same parameter on the boundary of the Mandelbrot set. 
To be precise, one starts by defining a relation on $\mathbb{Q}/\mathbb{Z}$ by setting that $\theta_1 \sim_\mathbb{Q} \theta_2$ if $R_M(\theta_1)$ 
and $R_M(\theta_2)$ land at the same point. Then, one takes the transitive closure of this relation and finally its topological closure 
to define an equivalent relation $\sim$ on $\mathbb{R}/\mathbb{Z}$. As constructed by Thurston \cite{Th}, there is a lamination $QML$ 
on the unit disk such that its induced equivalent relation is precisely $\sim$. If MLC holds, then $\theta_1 \sim \theta_2$ if and only it 
$R_M(\theta_1)$ and $R_M(\theta_2)$ land on the same point. 

Ray pairs are partially ordered: in fact, we say $(\theta_1^-, \theta_1^+) \prec (\theta_2^-, \theta_2^+)$ 
if the leaf $(\theta_1^-, \theta_1^+)$ separates $(\theta_2^-, \theta_2^+)$ from $0$. 
Let us denote as $N(\theta)$ the number of ends of the Hubbard tree associated to the angle $\theta$. 
Recall that if $\theta_1 \prec \theta_2$, then $N(\theta_1) \leq N(\theta_2)$. 

\medskip
A dyadic number $\theta_0$ defines a combinatorial vein in the Mandelbrot set as follows.

\begin{definition}
Given a dyadic rational number $\theta_0$, we define the \emph{combinatorial vein} of $\theta_0$ as the set of ray pairs $(\theta_1, \theta_2)$ 
such that: 
\begin{enumerate}
\item $\theta_0$ is the pseudocenter of $(\theta_1, \theta_2)$; 
\item $N(\theta_0) = N(\theta_1)$.
\end{enumerate}
\end{definition}

For instance, if $\theta_0 = \frac{1}{4}$, then the vein extends all the way to $(\frac{1}{7}, \frac{2}{7})$, 
since $N(\frac{1}{4}) = 3$, and one can check easily that $N(\frac{1}{7}) = 3$ (the ``rabbit"). 
On the other hand, if $\theta_0 = \frac{3}{16}$, then the vein extends up to $(\frac{39}{224}, \frac{43}{224})$. 
Indeed, $N(\frac{3}{16}) = 5$, and $N(\theta) = 3$ if $\theta \in (\frac{1}{7}, \frac{39}{224})$. 

\begin{figure} \label{F:not-simple}
\includegraphics[width = 0.7 \textwidth]{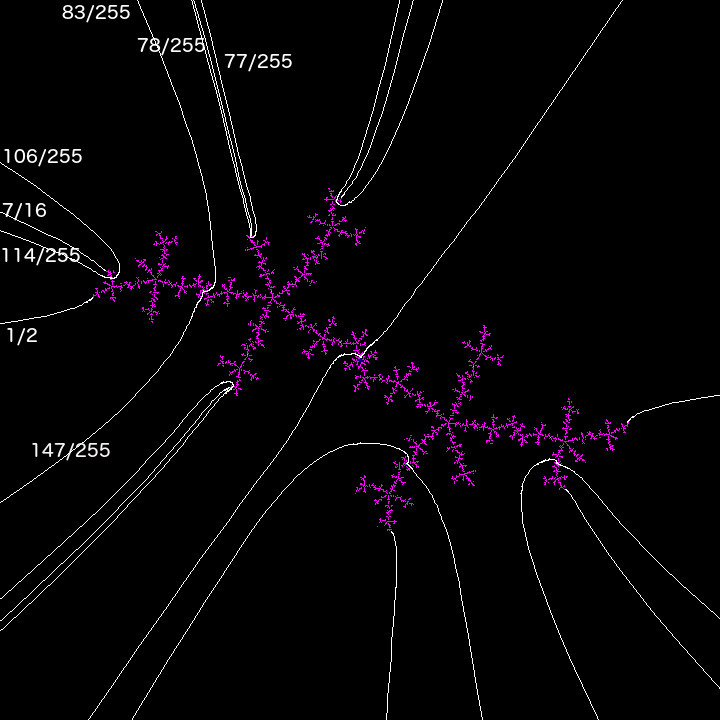}
\caption{The Julia set for the ray pair $(\theta^-, \theta^+) = (\frac{77}{255}, \frac{78}{255})$. One can see that the Hubbard tree for $\theta^-$ and $\theta^+$ has $6$ ends, while the Hubbard tree for the pseudocenter $\theta_0 = \frac{39}{128}$ has $8$ ends. Here, $\Vert \theta_0 \Vert = 7$, so the vein $V$ associated to $\theta_0$ has $\delta_V = 6$, but $\theta^+$ does not belong to $V$. This shows that the formula of Proposition \ref{P:index} cannot be extended beyond the vein $V$: indeed, the point in the forward orbit of $\theta^+$ which is closest to $\frac{1}{2}$ is $\frac{114}{255} \in \widehat{\mathcal{R}}$, while $D^{\delta_V}(\theta^+) = \frac{147}{255}$, so the formula would break down.}
\end{figure}

Note that condition (2) in the above definition is needed. For instance, if one considers the ray pair $(\theta_1, \theta_2) = (\frac{77}{255}, \frac{78}{255})$, then 
its pseudocenter is $\theta_0 = 39/128$. However, one checks that $N(\frac{77}{255}) = 6 < 8 = N(\frac{39}{128})$ (see Figure \ref{F:not-simple}).

\begin{definition}
The \emph{complexity} of a vein $V$ with pseudocenter $\theta_0$ is given by 
$$\delta_V = \Vert \theta_0 \Vert - 1.$$ 
\end{definition}

The complexity $\delta_V$ equals: 
\begin{itemize}
\item the smallest $k$ such that $D^k(\theta_0) = 1/2$; 
\item the smallest $k$ such that $D^k(\theta^-)$ and $D^k(\theta^+)$ lie on opposite connected components 
of the set $\mathbb{T}\setminus \{0, 1/2\}$;
\item the smallest $k$ such that $f^k_c(c)$ lies on the spine $[-\beta, \beta]$ of the Julia set
of $f_c$ for $c$ which belongs to the vein.
\end{itemize}

Clearly, $\delta_V = 0$ if and only if $V$ is the real vein. Given a vein $V$ in the upper half plane (i.e. with $\theta_0 <1/2$), its \emph{lower side} is the set of angles $\theta^+$ for which there exists $\theta^- < \theta^+$ such that $(\theta^-, \theta^+)$ belongs to the vein.

\begin{proposition}\label{P:index}
Let $V$ be a vein in the upper half plane, and $\delta_V$ be its complexity. 
Then for each ray pair $(\theta^-, \theta^+)$ which belongs to $V$ we have 
$$D^{\delta_V}(\theta^+) \subseteq \widehat{\mathcal{R}}.$$
\end{proposition}

In order to prove the proposition, we need the following

\begin{lemma} \label{L:side}
Let $(\theta^-, \theta^+)$ be a ray pair, and $\theta_0$ be its pseudocenter, with $\theta_0 < 1/3$ 
(i.e. the hyperbolic component lies in the upper half plane, and not in the $1/2$ limb). 
Then we have 
$$\theta^+ - \theta_0 < \theta_0 - \theta^-.$$
\end{lemma}

\begin{proof}
Recall that, for each rational $\frac{p}{q}$ with $p, q$ coprime, there exists a unique set $C_{p/q}$ of $q$ points on $\mathbb{T}$ such that the 
doubling map acts on $C_{p/q}$ with rotation number $p/q$, i.e. $D(x_i) = x_{i+p}$, 
where the indices are considered modulo $q$.
We know that the set $C_{p/q}$ is precisely the set of external rays which land on the $\alpha$ fixed point for $f_c$, 
when $c$ belongs to the $p/q$-limb.
Let us call \emph{sectors} the $q$ connected components of $\mathbb{T} \setminus C_{p/q}$. In particular, we call 
\emph{critical sector} $\mathcal{S}_1$ the smallest sector, and \emph{central sector} $\mathcal{S}_0$ the largest sector. 
We say an interval $I \subseteq \mathbb{T}$ is \emph{embedded} if it is entirely contained within one sector. 
Let us consider the interval $I_0 = [\theta^-, \theta^+]$, which is embedded by construction. 

Let $k \geq 0$ be the smallest integer such that $I_{k+1} := D^{k+1}(I_0)$ is not embedded (it must exist, 
since $D$ doubles lengths, and the length of an embedded interval is bounded above by a constant $c < 1$).
Note that the doubling map is a homeomorphism between each non-central sector and its image, 
so $I_k$ must be contained in the central sector. 

Let us now consider the set  $\widehat{C}_{p/q} = C_{p/q} + 1/2$, 
which is a subset of the central sector and equals the set of external rays which land on $-\alpha$, the preimage of the $\alpha$ fixed point. 
Thus, let us call \emph{subsectors} the connected components of $\mathcal{S}_0 \setminus \widehat{C}_{p/q}$.
In particular, exactly two subsectors $\mathcal{S}^-$ and $\mathcal{S}^+$ map to the critical sector $\mathcal{S}_1$; let us call them 
\emph{central subsectors}. 
Since the rays $D^k(\theta^-)$ and $D^k(\theta^+)$ land together, then either they are both contained in the same subsector, 
or one is contained in one central subsector and the other one in the other. 
If they are contained in the same subsector, then the whole $I_k$ is contained in the same subsector, hence its image $I_{k+1}$ 
is still embedded, contradicting the definition of $k$. 
Thus, the two endpoints of $I_k$ are contained in two distinct central subsectors. Since $I_k$ is embedded, then it must contain $0$; 
moreover, since $0$ is the dyadic of lowest possible complexity, this means that $D^k(\theta_0) = 0$.
Now, by construction the images $D^{k+1}(\theta^-)$ and $D^{k+1}(\theta^+)$ lie in the critical sector $\mathcal{S}_1$, 
which is contained inside the arc $[0, 1/2]$ since the hyperbolic component lies in the upper half plane.
This means that 
$$\ell([0, D^{k+1}(\theta^+)]) < 1/2 < \ell([D^{k+1}(\theta^-), 0])$$
(where $\ell$ denotes the length of the intervals) and, since $D^{k+1}$ is a homeomorphism when restricted to $[\theta_0, \theta^+]$ and to $[\theta^-, \theta_0]$,
$$2^{k+1} \ell([\theta_0, \theta^+]) = \ell(D^{k+1}[\theta_0, \theta^+]) = \ell([0, D^{k+1}(\theta^+)]) < 1/2$$
and similarly 
$$2^{k+1} \ell([\theta^-, \theta_0]) = \ell(D^{k+1}[\theta^-, \theta_0]) = \ell([D^{k+1}(\theta^-), 0]) > 1/2$$
hence by comparing the previous two equations
$$\ell([\theta_0, \theta^+]) <  \ell([\theta^-, \theta_0]) $$
which proves the claim.
\end{proof}

\subsection{Combinatorial Hubbard trees}

Recall that every angle $\theta \in \mathbb{T}$ has an associated lamination on the disc which is invariant by the doubling map (see \cite{Th}). 
The (two) longest leaves of the lamination are called \emph{major leaves}, and their common image is called \emph{minor leaf} 
and will be denoted by $m$. Moreover, we let $\beta$ denote the leaf $\{0 \}$, which we will take as the root of the lamination (the notation is due to the fact 
that the ray at angle $0$ lands at the $\beta$-fixed point). The dynamics on the lamination is induced by the dynamics of the doubling map 
on the boundary circle. In particular, let us denote $f : \overline{\mathbb{D}} \to \overline{\mathbb{D}}$ to be a continuous function 
on the filled-in disc which extends the doubling map on the boundary $S^1$. Let us denote by $\Delta$ the diameter of the circle which 
connects the boundary points at angles $\theta/2$ and $(\theta+1)/2$. Then we shall also choose $f$ so that it maps homeomorphically 
each connected component of $\mathbb{D}\setminus \Delta$ onto $\mathbb{D}$.

Let $\mathcal{L}_1$ and $\mathcal{L}_2$ be two distinct leaves. Then we define the \emph{combinatorial segment} $[\mathcal{L}_1, \mathcal{L}_2]$
as the set of leaves $\mathcal{L}$ of the lamination which separate $\mathcal{L}_1$ and $\mathcal{L}_2$. Some simple properties of 
combinatorial segments are the following: 

\begin{itemize}
\item[(a)]
if $\mathcal{L} \in [\mathcal{L}_1, \mathcal{L}_2]$, then $[\mathcal{L}, \mathcal{L}_1] \subseteq [\mathcal{L}_1, \mathcal{L}_2]$;
\item[(b)]
for any choice of leaves $\mathcal{L}_1, \mathcal{L}_2, \mathcal{L}_3$ we have $[\mathcal{L}_1, \mathcal{L}_2] \subseteq [\mathcal{L}_3, \mathcal{L}_1] \cup [\mathcal{L}_3, \mathcal{L}_2]$;
\item[(c)]
the image of $[\mathcal{L}_1, \mathcal{L}_2]$ equals: 
$$f([\mathcal{L}_1, \mathcal{L}_2])  = 
\left\{ \begin{array}{ll} 
\lbrack f(\mathcal{L}_1), f(\mathcal{L}_2) \rbrack &  \textup{if }\Delta\textup{ does not separate }\mathcal{L}_1\textup{ and }\mathcal{L}_2, \\
\lbrack f(\mathcal{L}_1), m \rbrack \cup \lbrack f(\mathcal{L}_2), m \rbrack & \textup{if }\Delta\textup{ separates }\mathcal{L}_1\textup{ and }\mathcal{L}_2. 
\end{array}
\right.$$
\item[(d)]
in any case, for any leaves $\mathcal{L}_1, \mathcal{L}_1$ we have 
$$[f(\mathcal{L}_1), f(\mathcal{L}_2)]  \subseteq f([\mathcal{L}_1, \mathcal{L}_2]) \subseteq \lbrack f(\mathcal{L}_1), m \rbrack \cup \lbrack f(\mathcal{L}_2), m \rbrack.$$
\end{itemize}
Finally, we shall say that a set $S$ of leaves is \emph{combinatorially convex} if whenever $\mathcal{L}_1$ and $\mathcal{L}_2$ 
belong to $S$, then the whole set $[\mathcal{L}_1, \mathcal{L}_2]$ is contained in $S$.

We now define 
$$H_n := \bigcup_{0 \leq i \leq n} [\beta, f^i(m)]$$
and
$$H := \bigcup_{n \in \mathbb{N}} H_n.$$

We call the set $H$ the \emph{combinatorial Hubbard tree} of $f$, as it is a combinatorial version of the (extended) Hubbard tree.

\begin{lemma} \label{L:invariant}
The combinatorial Hubbard tree $H$ has the following properties.
\begin{enumerate}
\item The set $H$ is the smallest combinatorially convex set of leaves which contains $\beta$, $m$ and is forward invariant.
\item 
Let $N \geq 0$ be an integer such that $f^{N+1}(m) \in H_N$. Then we have 
$H = H_{N}$.
\end{enumerate}
\end{lemma}

Note that $N + 2$ coincides with the number of ends of the extended Hubbard tree.

\begin{proof}
(1) Let us check that $H$ is combinatorially convex. Suppose $\mathcal{L}_1, \mathcal{L}_2$ belong to $H$, so that say $\mathcal{L}_1 \in [\beta, f^i(m)]$ 
and $\mathcal{L}_2 \in [\beta, f^j(m)]$. Then 
$$[\mathcal{L}_1, \mathcal{L}_2] \subseteq [\mathcal{L}_1, \beta] \cup [\beta, \mathcal{L}_2] \subseteq [f^i(m), \beta] \cup [\beta, f^j(m)] \subseteq H$$
as required. In order to check that $H$ is forward invariant, note that by point (d) and using that $f(\beta) = \beta$ yields the inclusion
\begin{equation} \label{E:iterate}
f([\beta, f^i(m)]) \subseteq [f(\beta), m] \cup [m, f^{i+1}(m)] \subseteq [\beta, m] \cup  [\beta, f^{i+1}(m)]
\end{equation}
and the right-hand side is contained in $H$ by construction, so $f(H) \subseteq H$.

In order to check the minimality, let $S$ be a combinatorially convex, forward invariant set of leaves which contains $\beta$ and $m$. 
Then  by convexity $S$ contains $[\beta, m]$. 
Moreover, by forward invariance for any $i \geq 0$ we have 
$$S \supseteq f^i([\beta, m]) \supseteq [f^i(\beta), f^i(m)] = [\beta, f^i(m)]$$
thus $S \supseteq H$. 

(2) Note that by construction $H_N \subseteq H$, and the same proof as in (1) shows that $H_N$ is combinatorially convex. In order to prove the claim 
it is thus enough to check that $H_N$ is forward invariant, because then by minimality we get $H_N \supseteq H$ and the claim is proven.
To prove forward invariance, note that 
$$f(H_N) = \bigcup_{0 \leq i \leq N} f([\beta, f^i(m)]) \subseteq H_N \cup [\beta, f^{N+1}(m)]$$
and since $f^{N+1}(m) \in H_N$ then $[\beta, f^{N+1}(m)] \subseteq H_N$, proving the claim.
\end{proof}

\begin{lemma}
Let $m_1 < m_2$, and $H_1, H_2$ be the corresponding combinatorial Hubbard trees. Then 
$$H_1 \subseteq H_2.$$
As a corollary, 
$$\#Ends(T_1) \leq \#Ends(T_2).$$
\end{lemma}

\begin{proof}
By definition, $m_1 < m_2$ means that $m_1 \in [\beta, m_2]$. Thus, $m_1 \in H_2$ and since $H_2$ is forward invariant we have $f^i(m_1) \in H_2$ 
for any $i \geq 0$.  Since also $\beta \in H_2$ and $H_2$ is combinatorially convex, then $[\beta, f^i(m_1)] \subseteq H_2$ for any $i \geq 0$, thus 
$H_1 \subseteq H_2$ as required. For the corollary, note that since the trees are dual to the laminations, $T_1 \subseteq T_2$, 
and in general a connected subtree of a tree cannot have more ends than the ambient tree. 
\end{proof}

\begin{lemma} \label{L:disjoint-intervals}
Let $(\theta^-, \theta^+)$ be a ray pair, and $\theta_0$ its pseudocenter, with $\Vert \theta_0 \Vert = q$. Then 
$$\#Ends(T_{\theta^+}) = \#Ends(T_{\theta_0})$$
if and only if the arcs $I_k = (D^k(\theta^-), D^k(\theta^+))$ for $k = 0, \dots, q - 1$ are disjoint.
\end{lemma}

\begin{proof}
Note that the forward orbit of $\theta_0$ has cardinality $q+1$, and since its minor leaf is a point (and so are all its forward iterates), 
the number of ends of $T_{\theta_0}$ is also $q+1$.
Now, consider the minor leaf $m = (\theta^-, \theta^+ )$ and its forward iterates.
Note that by Lemma \ref{l:pseudoiterate} the arc $I_k$ for $k \leq q -1$ never contains $0$.
Thus,  the intervals $I_j$ and $I_k$ are disjoint if and only if neither $f^j(m) \in [\beta, f^k(m)]$ nor $f^k(m) \in [\beta, f^j(m)]$.
Now, if all the intervals are disjoint, then $f^k(m) \notin [\beta, f^i(m)]$ for any $0 \leq i < k \leq q-1$, hence the number of ends 
of $T_{\theta^+}$ is at least $q+1$, which implies it is exactly $q+1$ since by the previous Lemma it cannot exceed the number of ends of $T_{\theta_0}$.
Vice versa, if $0 \leq i < k \leq q- 1$ and $f^k(m) \notin [\beta, f^i(m)]$, then either $I_k$ is disjoint from $I_i$, or $f^i(m) \in [\beta, f^k(m)]$. 
This however is impossible, as it implies $I_k \subseteq I_i$ with $i < k$. 
Indeed, let us denote as $x_k$ and $x_i$ the pseudocenters of, respectively, $I_k$ and $I_i$. 
As the complexity of the pseudocenter decreases precisely by $1$ under iteration, since $i < k$ one has $\Vert x_i \Vert > \Vert x_k \Vert$.
However, since $I_k \subseteq I_i$, then $x_k \in I_i$, hence by definition of pseudocenter $\Vert x_k \Vert \geq \Vert x_i \Vert$, contradicting the previous
statement. 
\end{proof}

 \begin{lemma} \label{l:pseudoiterate}
Let $[\alpha, \beta]$ be an embedded arc in $S^1$ which does not contain $0$, and $\theta_0$ its pseudocenter, with $q = \Vert \theta_0 \Vert$.
Then for all $0 \leq k \leq q-1$, the arc $[D^k(\alpha), D^k(\beta)]$ is embedded and does not contain $0$. 
\end{lemma}

\begin{proof}
We claim that there exists a minimal $k$ such that $[D^k(\alpha), D^k(\beta)]$ contains $1/2$, and moreover all arcs $[D^h(\alpha), D^h(\beta)]$ for $0 \leq h \leq k$ are embedded. Indeed, the map $D$ doubles lengths of arcs of length less than $1/2$, hence eventually the length of one of its images 
is more than $1/2$. Let $k$ be minimal such that the length of $J = [D^k(\alpha), D^k(\beta)]$ is at least $1/2$. Note that all previous iterates 
are embedded arcs, moreover $J$ must contain either $0$ or $1/2$. However, by minimality of $k$, since the preimages of $0$ are $0$ and $1/2$ and the initial arc does not contain $0$, $J$ must contain $1/2$. Since $J$ contains $1/2$ and does not contain $0$, then its pseudocenter is $1/2$, 
and since $J = D^k([\alpha, \beta])$ one has $D^k(\theta_0) = 1/2$, hence $k = q-1$, which completes the proof of the claim.
\end{proof}

\begin{proof}[Proof of Proposition \ref{P:index}.]
Let $k = \delta_V$, and let $\theta^- < \theta^+$ be the two endpoints of the leaf. We need to show that $D^k(\theta^+)$ belongs to $\widehat{\mathcal{R}}$.
Let us set $I_0 := (\theta^-, \theta^+)$ and $I_h := D^h(I_0)$. Since the component belongs to $V$,  
we have by Lemma \ref{L:disjoint-intervals} that the intervals $I_h$ for $0 \leq h \leq k+1$ are pairwise disjoint. Moreover, by definition of $\delta_V$ the interval $I_k$ contains $1/2$.
Now if one looks at the lamination it follows that  $f^{k+1}(m) \in H_k$, so we are in the hypothesis of Lemma \ref{L:invariant}. 
Thus, all higher iterates $f^i(m)$ for $i \geq k+1$ belong to $H_k$. This means that the leaf $f^k(m)$ separates the point $\{1/2\}$ from all postcritical leaves $f^i(m)$ with $i \geq 0$. Thus, if we consider the orbit $\{ D^i(\theta^+), i \geq 0 \}$, we have that no iterate lies in the interior of $I_k$, 
and points in complement of $I_k$ are by Lemma \ref{L:side} at least at distance $|D^k(\theta^+) - 1/2|$ from $1/2$, hence 
the element in the orbit closest to $1/2$ is $D^k(\theta^+)$, so $D^k(\theta^+)$ belongs to $\widehat{\mathcal{R}}$, as required.
\end{proof}

\section{Renormalization and landing} \label{S:renorm}

We proved that any element in the image of $\Phi_H$ combinatorially lands on the real axis. To complete the proof of the main theorem, 
we need to show it actually lands. In order to do so, we will prove that it is not renormalizable, hence it lands by Yoccoz's theorem. 

\begin{proposition} \label{P:landing}
Let $H$ be a hyperbolic component which does not lie in the $1/2$-limb, and let $\theta \in \mathbb{R}/\mathbb{Z}$ be an irrational angle 
of an external ray which lands on the boundary of $H$. Then the external ray at angle $\Phi_H(\theta)$ is not renormalizable, 
hence the corresponding ray lands. 
\end{proposition}

The proof uses the concept of maximally diverse sequence from \cite{Sh}, which we discuss in the following section.

\subsection{Maximally diverse sequences}

\begin{definition}
A sequence $(s_n) \in \{0, 1\}^\mathbb{N}$ is \emph{maximally diverse} if the subsequences
$$(s_{i + n p})_{n \in \mathbb{N}}$$ 
with $p \geq 1$, $0 \leq i \leq p-1$ are all distinct. 
\end{definition}

\begin{lemma}
A sequence $(s_n)$ is maximally diverse if and only if the subsequences
$$(s_{i + np})_{n \in \mathbb{N}}$$ 
with $p \geq 1$, $i \geq 0$ are all distinct.
\end{lemma}

\begin{proof}
Suppose that the two subsequences $(s_{i + np})_{n \in \mathbb{N}}$ and $(s_{j + nq})_{n \in \mathbb{N}}$ with $p, q \geq 1$ and $i, j \geq 0$ are equal. 
Then for any $k \geq 1$ we also have that $(s_{i + n pk})_{n \in \mathbb{N}}$ and $(s_{j + n qk})_{n \in \mathbb{N}}$ are equal. Now, we can choose $k$ large enough so that 
$i < pk$ and $j < qk$; then $s$ is not maximally diverse by definition. 
\end{proof}

\begin{corollary} \label{L:subsequence}
If $s$ is maximally diverse, then for any $i \geq 0$ and $p \geq 1$ the subsequence
$$(s_{i + np})_{n \in \mathbb{N}}$$ 
is maximally diverse. 
\end{corollary}

\begin{lemma} \label{L:shift}
Let $s \in \{0, 1\}^\mathbb{N}$ and $s' = \sigma(s)$ its shift, i.e. $(s')_n = s_{n+1}$ for any $n$. 
Then $s$ is maximally diverse if and only if $s'$ is maximally diverse. 
\end{lemma}

\begin{proof}
If $s$ is not maximally diverse, then there exist $i, j \geq 0$ and $p, q \geq 1$ such that 
$s_{i+np} = s_{j + nq}$ for all $n \geq 0$. Hence, also $s_{i+p+np} = s_{j +q+ nq}$ for all $n \geq 0$.
Since $i + p \geq 1$ and $j+q \geq 1$, this also implies
$s'_{i+p-1+np} = s'_{j +q-1+ nq}$ for all $n \geq 0$, hence $s'$ is not maximally diverse. Vice versa, if 
$s'$ is not maximally diverse then there exist $i, j \geq 0$ and $p, q \geq 1$ such that 
$s'_{i+np} = s'_{j + nq}$ for all $n \geq 0$, which implies 
$s_{i+1+np} = s_{j+1 + nq}$ for all $n \geq 0$, hence $s$ is not maximally diverse by the previous Lemma.
\end{proof}

Recall an infinite sequence $(\epsilon_n)$ is \emph{Sturmian} if there exists $\alpha \in (0, 1) \setminus \mathbb{Q}$, $\beta \in \mathbb{R}$ such that 
$$\epsilon_n = \lfloor (n+1) \alpha + \beta \rfloor - \lfloor n \alpha + \beta \rfloor - \lfloor \beta \rfloor$$
for all $n \geq 0$.

\begin{theorem}[Shallit \cite{Sh}] \label{T:Shallit}
Sturmian sequences are maximally diverse.
\end{theorem}

\subsection{Proof of Proposition \ref{P:landing}}

Let $\theta$ be an external angle of a ray landing on the boundary of the hyperbolic component $H$, of period $p$. 
If $\theta$ is irrational, then $\theta$ is the tuning of an irrational angle $\eta$ of a ray landing on the main cardioid. 
Thus, the binary expansion of $\theta$ is 
$$\theta := S_{\epsilon_0} S_{\epsilon_1} S_{\epsilon_2} \dots$$
where $(S_0, S_1)$ are the binary expansions of the two rays landing at the root of $H$, 
and $(\epsilon_n)$ is the binary expansion of the angle of a ray landing on the main cardioid. 
Thus, by \cite{BS} the sequence $(\epsilon_n)$ is a Sturmian sequence, hence by Theorem \ref{T:Shallit} it is a maximally diverse sequence. 

Hence, for any $\alpha = 0, \dots, p-1$
we have that either $(\theta_{\alpha + np})_{n \in \mathbb{N}}$ is maximally diverse (if $(S_0)_\alpha \neq (S_1)_\alpha$)
or is constant (if $(S_0)_\alpha = (S_1)_\alpha$). 
Thus, if we consider the image $\sigma := \Phi_H(\theta)$ we have by Lemma \ref{L:shift} that either 
$(\sigma_{\alpha + np})_{n \in \mathbb{N}}$ is maximally diverse 
or is eventually constant. 

Note that since $H$ does not intersect the real axis, then either $(S_0)_1 = (S_1)_1 = 0$ (if $H$ lies in the upper half plane) 
or $(S_0)_1 = (S_1)_1 = 1$ (if $H$ lies in the lower half plane). In both cases, the subsequence 
$(\theta_{1 + np})_{n \in \mathbb{N}}$
is eventually constant. Thus, there exists $j \geq 1$ such that the subsequence
$$(\sigma_{j + np})_{n \in \mathbb{N}}$$
is eventually constant. 

Let us now suppose by contradiction that the angle $\sigma$ is renormalizable. This implies that there exist two words $Z_0, Z_1$ 
with $|Z_0| = |Z_1| = q$ such that $\sigma \in \{ Z_0, Z_1 \}^\mathbb{N}$. Since by construction the angle $\sigma$ 
is real, the only possibility is that $(Z_1)_i = 1 - (Z_0)_i$ for any $i = 1, \dots, q$. 

\medskip

\textbf{Case 1.} Let us first assume that $p$ is not a multiple of $q$. 

\medskip

Then let $l := \textup{lcm}(p, q)$, and consider all remainder classes modulo $l$. 

\begin{definition}
Let us define two remainder classes $\alpha, \alpha'$ modulo $p$ to be $q$-equivalent if there exists $m \in \{0, \dots, \frac{l}{q}-1\}$ 
and $\beta, \beta' \in \{0, \dots, q-1\}$ such that $\alpha \equiv qm + \beta \mod p$ and $\alpha' \equiv qm + \beta' \mod p$. 
\end{definition}

\begin{lemma} \label{L:ec}
Let $p, q \geq 1$ be two integers, with $p$ not a multiple of $q$. 
Now, suppose that a set  $A \subseteq \{ 0, \dots, p-1 \}$ is not empty and has the following property: 
if $\alpha \in A$ and $\alpha' \in A$ is $q$-equivalent to $\alpha$, then $\alpha' \in A$. 
Then $A = \{ 0, \dots, p-1\}$. 
\end{lemma}

\begin{proof}
If $p \leq q$, the claim is almost trivial: indeed, the set $\{0, \dots, q-1\}$ projects to all possible classes modulo $p$. 
Let us suppose now $p > q$, and let $k := \lfloor \frac{p}{q} \rfloor$. Then each interval $A_m := [qm, qm+ q-1]$ with $0 \leq m \leq k$
lies in some equivalence class. Moreover, if we choose $j \in \{0, \dots, p-1\}$ such that $j \equiv (k+1)q \mod p$, then also 
each interval $B_m := [qm + j, qm + j + q-1]$ with $0 \leq m \leq k$ lies in some equivalence class. Then, note that 
$B_m$ intersects both $A_m$ and $A_{m+1}$, hence all remainder classes must belong to the same equivalence class.
\end{proof}

Now, let us pick $\alpha \in \{0, \dots, p-1\}$ such that $(\sigma_{\alpha + np})$ is eventually constant. 
Find $\gamma$ in $\{0, \dots, l-1\}$ such that $\gamma \equiv \alpha \mod p$, and let $\beta \in \{0, \dots, q-1\}$
such that $\beta \equiv \gamma \mod q$. Then one can write $\gamma = \beta + m q$ with $0 \leq m \leq \frac{l}{q} - 1$. 

Suppose that one has 
$$\sigma = Z_{\epsilon_0} Z_{\epsilon_1} \dots$$
with $\epsilon_i \in \{0, 1\}$. 

If for some $\alpha \in \{ 0, \dots, p-1 \}$ the sequence $(\sigma_{\alpha + np})_{n \in \mathbb{N}}$ is eventually constant (e.c.), so is its subsequence $(\sigma_{\gamma + n l})_{n \in \mathbb{N}}$,
which coincides with 
$((Z_{\epsilon_{m + nl/q}})_\beta)_{n \in \mathbb{N}}$. Since $(Z_0)_\beta \neq (Z_1)_\beta$, this implies 
$(\epsilon_{m + nl/q})_{n \in \mathbb{N}}$ is also eventually constant. On the other hand, if $(\sigma_{\alpha + np})_{n \in \mathbb{N}}$ is not eventually constant, then it is maximally diverse, 
hence by Lemma \ref{L:subsequence} so is 
$(\sigma_{\gamma + n l})_{n \in \mathbb{N}}$, hence also $(\epsilon_{m + nl/q})_{n \in \mathbb{N}}$ is not eventually constant.

Let us now consider another $\alpha' \in \{0, \dots, p-1 \}$ which is $q$-equivalent to $\alpha$. Then by the above discussion
$$(\sigma_{\alpha + np})_{n \in \mathbb{N}} \textup{ e.c. } \Leftrightarrow (\epsilon_{m + nl/q})_{n \in \mathbb{N}} \textup{ e.c. } \Leftrightarrow 
(\sigma_{\alpha' + np})_{n \in \mathbb{N}} \textup{ e.c. }$$

By Lemma \ref{L:ec}, since we know that there exists at least some $\alpha$ for which $(\sigma_{\alpha + np})_{n \in \mathbb{N}}$ is e.c., 
then each sequence $(\sigma_{\alpha'+np})_{n \in \mathbb{N}}$ is e.c. for any $\alpha' \in \{0, \dots, p-1\}$, hence $\sigma$ is also eventually periodic, 
thus the angle $\theta$ cannot be irrational. 

\medskip

\textbf{Case 2.} Finally, let us consider the case when $p$ is a multiple of $q$. 

\medskip

Then recall one can write 
$$\sigma = P S_{\epsilon_0} S_{\epsilon_1} \dots = Z_{\eta_0} Z_{\eta_1} \dots$$
where $P$ is a finite word of some length $k \geq 0$. 

Now, let us suppose that $k$ is not a multiple of $q$. 
Note that since $H$ is contained in either the upper or the lower half plane, we have $(S_0)_1 = (S_1)_1$; moreover, we claim that 
$(S_0)_p \neq (S_1)_p$: in fact, if one considers the external angles $\theta_0 = .\overline{S_0}$ and $\theta_1 = .\overline{S_1}$, 
the number of periodic external rays of period which divides $p$ and lie in the interval $(\theta_0, \theta_1)$ is even, 
since on each landing point exactly two rays land.  Hence, in their binary expansion the last digit of $S_0$ must be opposite to the last digit of $S_1$. 

Now, for each $i \geq 0$ there is some index $j$ such that $Z_{\eta_j}$ overlaps with both $S_{\epsilon_i}$ and $S_{\epsilon_{i+1}}$. 
As the last part of $Z_i$ must coincide with the first part of $S_{\epsilon_{i+1}}$, and since $S_0$ and $S_1$ start with the same symbol, 
this forces $Z_{\eta_j}$ to be either $Z_0$ or $Z_1$, independently of $i$. However, since the last symbols of $S_0$ and $S_1$ are different, this means that $S_{\epsilon_i}$ is also fixed, hence the sequence $(\epsilon_i)$ must be eventually constant, which contradicts the irrationality of $\theta$. 

Finally, if $k$ is multiple of $q$, then $S_0$ and $S_1$ are finite concatenations of $Z_0, Z_1$, which means that $\theta$ lies already 
in the small copy of the Mandelbrot set with roots $\{Z_0, Z_1\}$. Since $\{Z_0, Z_1\}$ represent a real pair, such a small copy of the Mandelbrot set lies in the $1/2$-limb, which contradicts the fact that $H$ is outside such limb and completes the proof of Proposition \ref{P:landing}. 

\subsection{Proof of the magic formula}

\begin{proof}[Proof of Theorem \ref{T:main}]
If $\theta$ belongs to $\Xi_H$, then by Proposition \ref{P:BC-formula} the angle $\bold{B}_H \bold{A}_H \cdot \theta$ belongs to the upper part of the combinatorial vein $V$ on which $H$ lies.
Hence, by Proposition \ref{P:index} the angle  $\Phi_H(\theta) := D^{\delta_V}( \bold{B}_H \bold{A}_H \cdot \theta)$ belongs to $\widehat{\mathcal{R}}$. Finally, 
by Proposition \ref{P:landing} the ray actually lands, hence $\Phi_H(\theta)$ belongs to $\mathcal{R}$, as claimed.
\end{proof}

\section{An alternate formula} \label{S:alter}

We conclude with another possible generalization of Douady's formula. This version does not depend on the vein structure of $\mathcal{M}$. 

\begin{proposition} \label{T:trivial-magic}
Let $H$ be any hyperbolic component, and $\Theta_H$ be the set of external angles of rays landing on the small Mandelbrot set with root $H$ (in particular, 
this set contains the set of external angles landing on the boundary of $H$). 
Then there exists an affine map $\varphi_H$ such that 
$$\varphi_H(\Theta_H) \subseteq \mathcal{R}.$$
\end{proposition}

\begin{lemma}
For each hyperbolic component $H$ of period $p > 1$, one has  
$$ \left|D^n (\theta) -\frac{1}{2} \right| \geq \frac{1}{2^{2p}}$$
for all $\theta \in \Theta_H$, for all $n \geq 0$.
\end{lemma}

\begin{proof}
Let $\Sigma_ 0 = \bold{A}_H$ and $\Sigma_1 = \bold{B}_H$ be the binary words which give the binary expansion of the two angles landing at the root of $H$, 
and let $p$ be the period of $H$, which equals the length of both $\Sigma_0$ and $\Sigma_1$.
By the construction of tuning operators, any external angle $\theta$ landing on the small Mandelbrot copy associated to $H$ has binary expansion of type
$$\theta = .\Sigma_{\epsilon_1} \Sigma_{\epsilon_2} \dots \Sigma_{\epsilon_n} \dots$$
where $\epsilon_i \in \{0,1\}$ for all $i$. Since $H$ is not the main cardioid, then both $\Sigma_0$ and $\Sigma_1$ contain both the symbol $0$ and the symbol $1$. As a consequence, any block of consecutive equal digits in the binary expansion of $\theta$ cannot have length larger than $2p - 2$.
However, all numbers in the interval $U_p = \left( \frac{1}{2} - \frac{1}{2^{2p}}, \frac{1}{2} + \frac{1}{2^{2p}}\right)$ have binary expansion of type either 
$$.0\underbrace{1\dots1}_{2p-1} \qquad \textup{or }\qquad .1\underbrace{0\dots0}_{2p-1},$$
hence none of the iterates $D^n(\theta)$ can lie in the interval $U_p$.
\end{proof}

\begin{proof}[Proof of Proposition~\ref{T:trivial-magic}]
The map is given by 
$$\varphi_H(\theta) := 0\underbrace{1\dots1}_{2p-1} \cdot \theta$$
In fact, by the above observation, $\varphi_H(\theta) \in U_p =  \left( \frac{1}{2} - \frac{1}{2^{2p}}, \frac{1}{2} + \frac{1}{2^{2p}}\right)$. On the other hand, consider the other iterates $D^n(\varphi_H(\theta))$ for $n \geq 1$.
If $n < 2p-1$, then $D^n(\varphi_H(\theta))$ begins with $11$ so it does not lie in $U_p$. For $n = 2p-1$, then $D^n(\varphi_H(\theta)) = 1 \cdot \theta$
also does not belong to $U_p$, as $\theta$ cannot begin with $2p-1$ zeros by the above Lemma. Finally, if $n \geq 2p$, then 
$$D^n(\varphi_H(\theta)) = D^{n-2p}(\theta) \notin U_p$$
again by the Lemma. In conclusion, since $\varphi_H(\theta)$ belongs to $U_p$ and none of its forward iterates does, then $\varphi_H(\theta)$ is closer 
to $1/2$ than all its iterates, hence $\varphi_H(\theta)$ belongs to $\mathcal{R}$.
\end{proof}

\end{document}